\documentclass[smallextended]{svjour3}       % onecolumn (second format)
\smartqed  % flush right qed marks, e.g. at end of proof
\usepackage{graphicx}
\usepackage{amssymb}
\usepackage{amsmath}
\usepackage[usenames]{color}
\usepackage{algorithm}
\usepackage[misc]{ifsym}
  {\end{list}}
\makeatletter

\gdef\listctr{list\romannumeral\the\@listdepth}\expandafter

\makeatother

  {\end{list}}

\newcommand{\R}{\mathbb{R}}
\newcommand{\N}{\mathbb{N}}

\newcommand{\ve}[1]{\boldsymbol{#1}}

\def\xk{\ve x^{(k)}}
\def\x{{\ve x}}
\def\y{{\ve y}}
\def\z{{\ve z}}
\def\p{{\ve p}}

\def\s{\ve\sigma} \def\sk{\ve \sigma^{(k)}}
\def\vesi{{\s}}%\sigma}
\def\vesik{\sk}%\vesi^{(k)}}
\def\hs{h_{\vesi}}
\def\hsk{h_{\vesik}}
\newcommand{\hsxy}[2]{\hs(#1,#2)}
\newcommand{\hskxy}[2]{\hsk(#1,#2)}

\newcommand{\hi}[2]{h^{#1}_{#2}}

\newcommand{\hixys}[4]{h^{#1}_{#2}(#3,#4)}
\def\ds{d_{\vesi}}

\newcommand{\dsxy}[2]{\ds(#1,#2)}

\def\xs{\ve {x^*}}
\def\barx{\bar{\ve x}}

\def\d{\ve d}
\def\arg{{\rm arg}}

\def\yk{\ve y^{(k)}}

\def\lamk{\lambda^{(k)}}

\def\xkk{\ve x^{(k+1)}}

\def\Hf{{\mathcal H}(f,\Omega,S)}

\def\endproof{\smartqed}%{\hfill$\square$\vspace{0.2cm}\\\noindent}

\def\endproof{\qed}%{\hfill$\square$\vspace{0.2cm}\\\noindent}

\journalname{}

\begin{document}

\title{A cyclic block coordinate descent method with generalized gradient projections}

%\titlerunning{Generalized gradient projection methods}

\author{Silvia Bonettini \and Marco Prato \and Simone Rebegoldi}

\institute{S. Bonettini \at
              Dipartimento di Matematica e Informatica, Universit\`a di Ferrara, Via Saragat 1, 44122 Ferrara, Italy \\
              Tel.: +39-0532-974785\\
              \email{silvia.bonettini@unife.it}
           \and
           M. Prato and S. Rebegoldi \at
              Dipartimento di Scienze Fisiche, Informatiche e Matematiche, Universit\`a di Modena e Reggio Emilia, Via Campi 213/b, 41125 Modena, Italy \\
              Tel.: +39-059-2055590\\
              \email{marco.prato@unimore.it, simone.rebegoldi@unimore.it}
}

\date{Received: date / Accepted: date}

\maketitle

\begin{abstract}
The aim of this paper is to present the convergence analysis of a very general class of gradient projection methods for smooth, constrained, possibly nonconvex, optimization. The key features of these methods are the Armijo linesearch along a suitable descent direction and the non Euclidean metric employed to compute the gradient projection. We develop a very general framework from the point of view of block--coordinate descent methods, which are useful when the constraints are separable.
\end{abstract}

\section{Introduction}\label{sec:intro}
This paper deals with the problem
\begin{equation}\label{minf}
\min_{\ve x\in \Omega} f(\ve x),
\end{equation}
where $\Omega\subseteq\R^n$ is a closed and convex set and $f$ is a continuously differentiable function. The aim of this work is to generalize the class of gradient projection methods whose basic iteration is given by
\begin{equation}\label{iter0}
\xkk = \xk + \lamk (\yk-\xk),
\end{equation}
where $\yk$ is the Euclidean projection of $\xk-\sigma_k\nabla f(\xk)$ onto $\Omega$, i.e.
\begin{equation}\label{y1}
\yk = P_\Omega(\xk-\sigma_k\nabla f(\xk) ) \equiv \arg\min_{\x\in\Omega}\|\x - \xk +\sigma_k\nabla f(\xk) \|^2,
\end{equation}
and $\sigma_k > 0$, $\lamk\in(0,1]$ control the steplength.\\
Iteration \eqref{iter0}--\eqref{y1} is also referred as gradient projection method \emph{with linesearch along the descent direction} \cite{Bertsekas,Iusem03}, which depends on two parameters $\lamk,\sigma_k$. Usually, in iteration \eqref{iter0}, $\lamk$ is adaptively computed to ensure the sufficient decrease of the objective function and, thus, the convergence of the whole scheme, while $\sigma_k$ is a 'free' parameter which can be chosen in order to improve the effectiveness of the algorithm (see e.g. \cite{Barzilai1988,Dai2006,DeAsmundis2013,Fletcher2012}).\\
In our analysis, we extend the convergence results about the gradient projection method \eqref{iter0}-\eqref{y1} to the more general case where $\yk$ is defined as
\begin{equation}\label{yh}
\yk =  \arg\min_{\x\in\Omega}\hskxy\x\xk
\end{equation}
and $\hs$ is a suitable strictly convex function depending on the array of parameters $\vesi\in \R^q$. The choice of $\hs$ can be addressed by taking into account some features of problem \eqref{minf}. For example, $\hskxy\x\xk$ may represent a local approximation of $f$ at $\xk$, or may play the role of barrier for a given constraint set, forcing the iterates to stay in the interior of it \cite{Auslender-etal-2007,Auslender-Teboulle-2005,Auslender-Teboulle-2006}.\\
In particular, we present our results in the more general framework of the block--coordinate methods, which are useful
when the constraint set in \eqref{minf} has a separable structure, i.e. $\Omega = \Omega_1\times\dots\Omega_m$, with $\Omega_i\subseteq \R^{n_i}$, $\sum_{i=1}^m n_i = n$, so that any $\x\in\Omega$ can be block partitioned as $\x = (\x_1^T,\dots,\x_m^T)^T$, $\x_i\in \R^{n_i}$.\\
Such methods are based on the idea of performing successive minimizations over each block, as in the classical nonlinear Gauss-Seidel method \cite{Bertsekas}:
\begin{equation}
\xkk_i\in \mbox{arg}\min_{\x\in\Omega_i}
f(\xkk_1,...,\xkk_{i-1},\x,\xk_{i+1},...,\xk_m).
\label{GS}
\end{equation}
However, the convergence of this approach is not ensured without quite restrictive convexity assumptions (see \cite{Bertsekas,Grippo00}) and, in addition, computing an exact minimum of $f$, even if restricted to a single block, can be impractical.\\
On the other side, inspired by the idea of \eqref{GS}, effective methods able to handle general nonconvex problems and with global convergence properties can be designed \cite{Bonettini11,Cassioli-etal-2013,Grippo1999}.\\
In this paper we further develop the cyclic block gradient projection method proposed in \cite{Bonettini11}, allowing generalized projections based on non Euclidean distances. In particular, we propose a method consisting in applying a finite number of iterations of the form \eqref{iter0}--\eqref{yh} to each subproblem of type \eqref{GS} and we show that any limit point of the generated sequence is stationary without any convexity assumption.
Our general framework includes, but it is not limited to, several state-of-the-art methods, such as the scaled gradient projection method \cite{BZZ09}, the spectral projected gradient method \cite{BMR03}, the cyclic block gradient projection method \cite{Bonettini11} and the successive convex approximation algorithm \cite{Razaviyayn-etal2013}.\\
The paper is organized as follows: in section \ref{sec:projections} we devise the property of the operator $\hs$ in \eqref{yh} which allow to reformulate the stationarity condition of \eqref{minf} by means of a class of generalized projection operators. We also show that they can be used to design families of descent directions. Building on this material and on the well known properties of the Armijo linesearch procedure, in section \ref{sec:blocks} we define a block--coordinate generalized gradient projection method and we develop the related convergence analysis. Our conclusions are given in section \ref{sec:concl}.
\section{Generalized gradient projections}\label{sec:projections}
In this section we give the definition of a generalized projection operator, providing some examples of well-known functions belonging to this category.
\begin{definition}\label{Def_projgrad}
Let $S\subseteq\R^q$. We define a metric function associated to $f$ any continuously differentiable function $\hs:\Omega\times\Omega\rightarrow \R$ such that for any choice of the parameter $\s \in S$ the following properties hold:
\begin{itemize}
\item[(H1)] $\hs$ is convex with respect to its first argument, i.e.
\begin{equation}\label{h_convex}
\hsxy{\y}{\z}\geq \hsxy{\x}{\z} + \nabla_1 \hsxy{\x}{\z}^T(\y-\x) \quad \forall \x,\y,\z\in\Omega
\end{equation}
and, for any $\z\in \Omega $, $\hsxy{\cdot}{\z}$ admits a unique minimum point;
\item[(H2)] for any point $\x\in\Omega$ and for any feasible direction $\ve d\in \R^n$ we have
\begin{equation}\label{H2}
\nabla_1\hsxy{\x}{\x}^T\ve d = \nabla f(\x)^T\d;
\end{equation}
\item[(H3)] $\hs$ continuously depends on the parameter $\s$.
\end{itemize}
We denote by $\Hf$ the set of the metric functions satisfying properties (H1)--(H3) and, for any $\hs\in\Hf$, we define the associated generalized gradient projection operator $\p(\ \cdot \ ;\hs):\Omega \rightarrow\Omega$ as
\begin{equation}\label{def:gen_proj_op}
\p(\x;\hs) = \arg\min_{\z\in \Omega} \hsxy{\z}{\x} \quad \forall\x\in \Omega.
\end{equation}
\end{definition}
\begin{example}\label{esempio1}
Properties \eqref{h_convex}--\eqref{H2} are satisfied when the function $\hs$ is defined as
\begin{equation}\label{h-d}
\hsxy{\x}{\y} = \nabla f(\y)^T(\x-\y) + \dsxy{\x}{\y},
\end{equation}
where $\ds\in\mathcal{D}(\Omega)$. In these settings we can find:
\begin{itemize}
\item[a)] the standard Euclidean projection $\p(\x;h_\sigma) = P_{\Omega}(\x-\sigma\nabla f(\x))$, obtained by choosing
\begin{equation} \label{eucl}
d_\sigma(\x,\y) = \frac 1 {2\sigma}\|\x-\y\|^2, \quad \sigma > 0;
\end{equation}
\item[b)] the scaled Euclidean projection, considered for example in \cite{BMR03,BZZ09}, corresponding to the choice
\begin{equation}\label{scaled_euclidean_projection}
d_{(\alpha,D)}(\x,\y) = \frac 1 {2\alpha}(\x-\y)^TD^{-1}(\x-\y).
\end{equation}
In this case the array of parameters $\s$ is given by the pair $(\alpha, D)$, where $\alpha\in \R_{>0}$ and $D\in \R^{n\times n}$ is a symmetric positive definite matrix;
\item[c)] the Bregman distance associated to a strictly convex function $b:\Omega\rightarrow\R$, which is defined as
\begin{equation}\label{bregman_function}
d_\sigma(\x,\y) = \frac{1}{\sigma} (b(\x) - b(\y) - \nabla b(\y)^T(\x-\y)), \quad \sigma > 0.
\end{equation}
\end{itemize}
\end{example}
\begin{example}\label{esempio2}
If $f$ is convex, a further class of functions satisfying the properties of Definition \ref{Def_projgrad} is given by
\begin{equation}\label{ex2_2}
\hsxy{\x}{\y} = f(\x) + \dsxy{\x}{\y},
\end{equation}
where again $\ds\in\mathcal{D}(\Omega)$. If $\ds$ is chosen as in \eqref{eucl}, the resulting $\hs$ leads to the so-called \emph{resolvent (or proximity) operator} associated to $f$ (see e.g. \cite{Chambolle2011,Combettes-Pesquet-09} or \cite{Eckstein1993} for the general case).
\end{example}
\begin{example}\label{esempio3}
Consider the case when $f=f_0+f_1$, where $f_0,f_1:\Omega \rightarrow \R$ and $f_0$ is convex. Then, the function defined as
\begin{equation}\label{ex2_3}
\hsxy{\x}{\y} = f_0(\x) + \dsxy{\x}{\y} + \nabla f_1(\y)^T(\x-\y) \ \ \forall \x,\y \in \Omega,
\end{equation}
with $\ds\in\mathcal{D}(\Omega)$, belongs to $\Hf$. If $\ds$ reduces to \eqref{eucl}, the corresponding projection operator is also known as the \emph{proximal gradient operator}, which is employed to define forward-backward splitting algorithms for convex optimization \cite{Combettes-Pesquet-09,Razaviyayn-etal2013}.
\end{example}
Observe that the metric functions defined in \eqref{ex2_2}-\eqref{ex2_3} are majorant of the objective function, that is $\hsxy\x\y\geq f(\x)$ for all $\x,\y\in \Omega$.
Further, any convex \emph{upper bound function} in the sense of \cite[Assumption 1]{Razaviyayn-etal2013} admitting a unique minimum point clearly satisfies the premises of Definition \ref{Def_projgrad}.\\
{\it Remark.} For sake of simplicity, in Definition \ref{Def_projgrad} we assume $\hs$ and $f$ to be smooth functions, but this could be relaxed, requiring only the existence of directional derivatives. Indeed, properties \eqref{h_convex} and \eqref{H2}, as well as the analysis carried out in the rest of this section, could be reformed in terms of directional derivatives.\\
In general, any function $\hs\in \Hf$ can be exploited to define a descent direction for problem \eqref{minf}, as stated in the following proposition.
\begin{proposition}\label{Prograd}
Let $\x\in \Omega$, $\s\in S\subseteq\R^q$, $\hs\in\Hf$ and
\begin{equation}\label{puppa}
\y = \p(\x;\hs).
\end{equation}
Then we have that
\begin{equation}\label{thesis_prograd}
\nabla f(\ve x)^T (\y-\x)\leq 0
\end{equation}
and the equality holds if and only if $\y=\x$.
\end{proposition} \noindent
{\it Proof. }
Inequality \eqref{h_convex} with $\z=\x$ yields
\begin{equation}\nonumber
\nabla_1\hsxy{\x}{\x}^T(\y-\x) \leq \hsxy{\y}{\x}-\hsxy{\x}{\x}\leq 0,
\end{equation}
where the rightmost inequality follows from \eqref{def:gen_proj_op} and, since the minimum point of $\hsxy{\cdot}{\x}$ is unique, the equality holds if and only if $\x=\y$. Then, the thesis follows recalling \eqref{H2}.
\endproof
In the following proposition, we show that the stationary points of \eqref{minf} can be characterized as fixed points of the generalized projection operator \eqref{def:gen_proj_op}.
\begin{proposition}\label{Proproj}
Let $S \subseteq \R^q$, $\vesi\in S$ and $\hs\in\Hf$. A point $\ve x\in \Omega$ is a stationary point for problem \eqref{minf} if and only if $\x=\p(\x;\hs)$.\\
\end{proposition}\noindent
{\it Proof. }Assume that for a point $\xs\in \Omega$ the following equality holds:
\begin{equation*}
\xs = \arg\min_{\x\in\Omega} \hsxy{\x}{\xs}.
\end{equation*}
Then, the stationarity of $\xs$ yields
\begin{equation*}
\nabla_1 \hsxy{\xs}{\xs}^T(\x-\xs) \geq 0 \quad \forall \x\in \Omega.
\end{equation*}
Since by assumption \eqref{H2} we have $\nabla_1 \hsxy{\xs}{\xs}^T(\x-\xs) = \nabla f(\xs)^T(\x-\xs)$, it follows that $\xs$ is a stationary point for problem \eqref{minf}.\\
Conversely, let $\xs\in\Omega$ be a stationary point of \eqref{minf} and define
\begin{equation*}
\barx= \arg\min_{\x\in\Omega} \hsxy{\x}{\xs}.
\end{equation*}
Assume by contradiction that $\xs\neq\barx$. Then, combining \eqref{h_convex} with $\x=\z=\xs$, $\y=\barx$ and \eqref{H2} we obtain
\begin{equation*}
\nabla f(\xs)^T(\barx-\xs) \leq \hsxy{\barx}{\xs}-\hsxy{\xs}{\xs} < 0,
\end{equation*}
where the last inequality follows from the fact that $\barx$ is the unique minimum point of $\hsxy{\cdot}{\xs}$ and $\xs\neq\barx$. This contradicts the stationarity assumption on $\xs$.
\endproof
\section{Cyclic block generalized gradient projection method}\label{sec:blocks}
In this section we consider problem \eqref{minf} where the constraint set has the following separable structure
\begin{equation}\label{Omegai}
\Omega = \Omega_1\times\dots\Omega_m, \ \Omega_i\subseteq \R^{n_i},\  \sum_{i=1}^m n_i = n\end{equation}
so that any $\x\in\Omega$ can be block partitioned as $\x = (\x_1^T,\dots,\x_m^T)^T$, $\x_i\in \R^{n_i}$.\\
The key ingredients of our approach are the sufficient decrease of the objective function enforced by a block version of the well known Armijo backtracking procedure and a suitable metric function $\hs\in \Hf$ defined so that is separable with respect to the partition in \eqref{Omegai}.\\
Then, we first recall in Algorithm \ref{Algo2} the block version of the Armijo linesearch method.
\begin{algorithm}[ht]\caption{Armijo linesearch algorithm}\label{Algo2}\noindent
Let $\{\ve z^{(k)}\}_{k \in \N}$ be a sequence of points in $\Omega$ and $\{\ve d_i^{(k)}\}_{k \in \N}$ a sequence of descent directions, for a given $i\in\{1,...,m\}$. Fix $\delta_i, \beta \in(0,1)$ and compute $\lambda_i^{(k)}$ as follows:
\begin{itemize}
\item[1.] Set $\lambda_i^{(k)}=1$;
\item[2.] \textsc{If}
\begin{equation}
f(\ve z_1^{(k)},...,\ve z_i^{(k)}+\lambda_i^{(k)}\ve d_i^{(k)},...,\ve z_m^{(k)})\leq
f(\ve z^{(k)})+\beta \lambda_i^{(k)}\nabla_i f(\ve z^{(k)})^T\ve d_i^{(k)}
\label{Armijo}
\end{equation}
\textsc{Then}  go to step 3.\\
\textsc{Else}  set $\lambda_i^{(k)}=\delta_i\lambda_i^{(k)}$ and go to step 2.
\item[3.]\textsc{End}
\end{itemize}
\end{algorithm}
In the following proposition we give conditions which guarantee that Algorithm \ref{Algo2} is well defined. Its proof can be derived from known results (see \cite{Bertsekas,Grippo00}).
\begin{proposition}\label{ProArm}
Let $\{\ve z^{(k)}\}_{k \in \N}$ be a sequence of points in $\Omega$. Assume that $\ve z^{(k)}$ converges to some $\bar{\ve z}$ and for $i\in\{1,...,m\}$ let
$\{\ve d_i^{(k)}\}_{k \in \N}$ be a sequence of feasible directions such that
\begin{itemize}
\item[(A1)] there exists a number $M>0$ such that $\|\ve d^{(k)}_i\|\leq M$ for all $k \in \N$;
\item[(A2)] we have $\nabla_i f(\ve z^{(k)})^T \ve d_i^{(k)} < 0$ for all $k \in \N$;
\item[(A3)] we have $\displaystyle\lim_{k\rightarrow\infty}f(\ve z^{(k)})- f(\ve z_1^{(k)},...,\ve z_i^{(k)}+
\lambda_i^{(k)}\ve d_i^{(k)},...,\ve z_m^{(k)})
=0$, where $\lambda_i^{(k)}$ is computed with Algorithm \ref{Algo2}.
\end{itemize}
Then, for each $k\in\N$ the LS procedure terminates in a finite number of steps and, furthermore,
$\lim_{k\rightarrow\infty} \nabla_i f(\ve z^{(k)})^T\ve d^{(k)}_i = 0$.
\end{proposition}
In order to formally introduce the method and perform its convergence analysis, we choose the metric function $\hs\in \Hf$, where $S=S_1\times...\times S_m$, $S_i\subset \R^{q_i}$, such that the parameter $\s$ can be partitioned as $\s = (\s_1,\dots,\s_m)$. Moreover, we define $\hs$ so that it is separable over the $m$ blocks with respect to its first variable, i.e.
\begin{equation}\label{hp_separable}
\hsxy\x\y = \sum_{i=1}^m \hixys i{\s_i}{\x_i}\y,
\end{equation}
where the functions $\hi i{\vesi_i}:\Omega_i\times\Omega \rightarrow \R$ satisfy the following conditions:
\begin{itemize}
\item[(BH1)] $\hi i{\vesi_i}$ is convex with respect to its first argument and admits a unique minimum point;
\item[(BH2)] for any point $\x\in\Omega$ and for any vector $\ve d\in \R^{n_i}$ such that $\x_i + \d \in \Omega_i$ we have
\begin{equation}\label{BH2}
\nabla_1\hixys i{\vesi_i}{\x_i}\x^T\ve d = \nabla_i f(\x)^T\d,
\end{equation}
where $\nabla_i f(\x)$ denotes the gradient of $f$ with respect to the $i$--th block of variables;
\item[(BH3)] $\hi i{\vesi_i}$ continuously depends on the parameter $\vesi_i\in \R^{q_i}$.
\end{itemize}
It is easy to see that the metric function $\hs$ defined in \eqref{hp_separable}, thanks to the assumptions (BH1)--(BH3), belongs to $\Hf$ and the associated generalized gradient projection can be also partitioned by blocks as
\begin{equation}\label{blocchi}
\p(\x;\hs) = \begin{pmatrix} \p_1(\x;\hi 1{\vesi_1})\\ \vdots \\ \p_m(\x;\hi m {\vesi_m})\end{pmatrix}, \ \ \mbox{ where }\ \
\p_i(\x;\hi i{\vesi_i}) = \mbox{arg}\min_{\z_i\in \Omega_i} \hixys i{\vesi_i}{\z_i}\x.
\end{equation}
\begin{lemma}\label{lem:blocchi}
Let $\x\in\Omega$ and $\s\in S\subseteq \R^q$. Then,
\begin{itemize}
\item[(i)] $\x$ is stationary for problem \eqref{minf} if and only if $\p_i(\x;h_{\vesi_i}^i) = \x_i$ $\forall i=1,\dots, m$;
\item[(ii)]  $\nabla_i f(\x)^T(\p_i(\x;\hi i{\vesi_i})-\x_i)\leq 0$ $\forall i=1,\dots, m$ and the equality holds if and only if $\x_i = \p_i(\x;\hi i{\vesi_i})$.
\end{itemize}
\end{lemma}
Part (i) of the previous Lemma directly follows from \eqref{blocchi} and from Proposition \ref{Proproj}, while part (ii) can be easily proved by employing the same arguments as in the proof of Proposition \ref{Prograd}.\\

\begin{algorithm}[ht]\caption{Cyclic Block Generalized Gradient Projection Method}\label{AlgoBlocchi}
Define a compact set $S$ and a metric function $\hs\in\Hf$ as in \eqref{hp_separable}. Choose $\beta, \delta\in (0,1)$.\\
Choose $\ve x^{(0)}\in \Omega$ and the upper bounds for the inner iterations numbers $L_1,\dots,L_m$.\\
{\textsc{For}} $k=0,1,2,...$
\begin{itemize}
\item[1] {Set} $\ve z(k,0)=\ve x^{(k)}$
\item[2]\textsc{For} $i=1,...,m$
\begin{itemize}
%\item[2.1] Set $\ve d_i^{(k)}=\nabla_i^Pf(\ve z(k,i-1))$;
\item[2.1] Set $\x_i^{(k,0)}= \x_i^{(k)}$
\item[2.2] Choose the inner iterations number $L_i^{(k)}\leq L_i$
\item[2.3] \textsc{For} $\ell = 0,...,L^{(k)}_i-1$
    \begin{itemize}
    \item[2.3.0] Set $\tilde \x^{(k,\ell)} = (\xkk_1,\dots,\xkk_{i-1},\x_i^{(k,\ell)},\xk_{i+1},\dots,\xk_m)$
    \item[2.3.1] Choose the parameter $\vesi_i^{(k,\ell)}\in S_i$
    \item[2.3.2] Compute the descent direction
    \begin{equation*}
		\d^{(k,\ell)}_i = \p_i(\tilde \x^{(k,\ell)} ; \hi i {\vesi_i^{(k,\ell)}}) - \x_i^{(k,\ell)}
		\end{equation*}
and set $\tilde \d^{(k,\ell)} = (0,\dots,0,\d_i^{(k,\ell)},0,\dots,0)$
    \item[2.3.3] Compute with Algorithm \ref{Algo2} the Armijo steplength $\lambda^{(k,\ell)}_i$ such that
        \begin{equation*}
        f(\tilde \x^{(k,\ell)}+\lambda^{(k,\ell)}_i\tilde \d^{(k,\ell)}) \leq f(\tilde \x^{(k,\ell)})+ \beta \lambda^{(k,\ell)}_i \nabla_i f(\tilde \x^{(k,\ell)})^T\d^{(k,\ell)}_i
        \end{equation*}
    \item[2.3.4] Set $\x_i^{(k,\ell+1)} =  \x_i^{(k,\ell)}+\lambda^{(k,\ell)}_i\ve d^{(k,\ell)}_i $
    \end{itemize}
    \textsc{End}
\item[2.4] Set $ \x_i^{(k+1)} = \x_i^{(k,L_i^{(k)})}$
\item[2.5] Set $\z(k,i)=(\xkk_1,...,\xkk_i,\xk_{i+1},...,\xk_m)$
\end{itemize}
\textsc{End}
\item[3] Set $\ve x^{(k+1)}=\ve z(k,m)$
\end{itemize}
\textsc{End}
\end{algorithm}

\noindent The previous results can be exploited to design a cyclic block generalized gradient projection (CBGGP) method, whose steps are outlined in Algorithm \ref{AlgoBlocchi}. Before to analyze the convergence properties of this approach, we observe that it is a descent method and, in particular, the objective function is nondecreasing over the {\it partial updates} $\z(k,i)$, $i=0,...,m$, $k=1,2,...$ defined at step 2.5. Indeed, the following inequalities hold
\begin{equation*}
f(\z(k,i+1))\leq f(\z(k,i))+\beta\lambda_{i+1}^{(k,0)}\nabla_{i+1}f(\z(k,i))^T\d_{i+1}^{(k,0)}\leq f(\z(k,i))
\end{equation*}
which also implies
\begin{align}\label{monotone_z}
f(\z(k+1,0)) = f(\z(k,m)) &\leq f(\z(k,i+1)) \\
                          &\leq f(\z(k,i)) \leq f(\z(k,0)) = f(\z(k-1,m)). \nonumber
\end{align}
We are now ready to give the first result about Algorithm \ref{AlgoBlocchi}.
\begin{proposition}\label{problocchi}
Let $\{\ve x^{(k)}\}_{k \in \N}$ be the sequence generated by Algorithm \ref{AlgoBlocchi}. Suppose that for some $i\in\{0,...,m\}$ the sequence $\{\ve z(k,i)\}_{k \in \N}$ admits a limit point $\ve{\bar z}$. Then $\p_{i+1}(\ve{\bar z}; h^{i+1}_{\s_{i+1}}) = \ve{\bar z}_{i+1}$ $\forall \s_{i+1}\in S_{i+1}$ if $i<m$, while $\p_{1}(\ve{\bar z}; h^{1}_{\s_{1}}) = \ve{\bar z}_{1}$ $\forall \s_{1}\in S_{1}$ if $i=m$.
\end{proposition} \noindent
{\it Proof.} Suppose first that $i<m$. From Lemma \ref{lem:blocchi}, we only need to show that there exists $\bar \s_{i+1}\in S_{i+1}$ such that equality $\p_{i+1}(\ve{\bar z}; h^{i+1}_{ \bar \s_{i+1}}) = \ve{\bar z}_{i+1}$ holds.\\
Assume by contradiction that $\p_{i+1}(\ve{\bar z}; \hi {i+1}{\vesi_{i+1}}) \neq \ve{\bar z}_{i+1}$ for all $\s_{i+1}\in S_{i+1}$.
Let $K$ be the set of indices such that $\{\ve z(k,i)\}_{k\in K}$ converges to $\bar{\ve z}$ and $\{\s_{i+1}^{(k,0)}\}_{k\in K}$ converges to some $\bar \s_{i+1}\in S_{i+1}$. If $\|\p_{i+1}(\ve{\bar z}; \hi{i+1}{\bar \vesi_{i+1}}) - \ve{\bar z}_{i+1}\|=2\epsilon >0$, the continuity of the generalized projection operator with respect to all its arguments guarantees that, for
$k\in K$ being sufficiently large, we have
\begin{equation*}
\|\d^{(k,0)}_{i+1}\| > \epsilon > 0,
\end{equation*}
where $\d^{(k,0)}_{i+1}=\p_{i+1}(\z(k,i);\hi{i+1}{\vesi_{i+1}^{(k,0)}})-\xk_{i+1}$  (see also Step 2.3.2 of Algorithm \ref{AlgoBlocchi}).
Then, by applying Lemma \ref{lem:blocchi} (ii) we have
\begin{equation}\label{contr}
\nabla_{i+1} f(\ve z(k,i))^T\ve d^{(k,0)}_{i+1}\leq -\eta<0,
\end{equation}
where $\eta$ is some positive scalar.\\
On the other side, inequalities \eqref{monotone_z} guarantee that, for all $i$, we have $\lim_{k\rightarrow \infty} f(\z(k,i))=f(\bar {\z})$, thus we obtain that
\begin{equation*}
\lim_{k\rightarrow \infty}f(\ve z(k,i))-f(\xkk_1,...,\xkk_i,\xk_{i+1}+ \lambda_{i+1}^{(k,0)}\d^{(k,0)}_{i+1},...,\xk_m )=0.
\end{equation*}
Moreover, since $\{\ve z(k,i)\}_{k\in K}$ is a convergent sequence, it is also bounded. Therefore the sequence $\{\ve d^{(k,0)}_{i+1}\}_{k\in K}$ is bounded and Proposition \ref{ProArm} implies that
\begin{equation*}
\lim_{k\rightarrow\infty,k\in K}\nabla_{i+1} f(\ve z(k,i))^T\ve d^{(k,0)}_{i+1}=0,
\end{equation*}
which contradicts \eqref{contr}.\\
The same arguments can be applied also when $i=m$,
since $\ve z(k,m)=\ve z(k+1,0)$.
\endproof
The previous proposition is crucial for proving the main convergence result for Algorithm \ref{AlgoBlocchi}, given below.
\begin{theorem}
Let $\{\ve x^{(k)}\}_{k \in \N}$ be the sequence generated by Algorithm \ref{AlgoBlocchi} and assume that $\bar{\ve x}$ is a limit point of $\{\ve x^{(k)}\}_{k \in \N}$. Then $\bar{\ve x}$ is a limit point also for the sequences $\{\ve z(k,i)\}_{k \in \N}$ for any $i=1,...,m-1$ and it is a stationary point for problem \eqref{minf}.
\end{theorem} \noindent
{\it Proof.} The proof runs by induction on the block index $i$ and on the inner iteration number $\ell$ and it is similar to that of Theorem 4.2 in \cite{Bonettini11}. Since $\bar{\ve x}$ is a limit point for $\{\ve x^{(k)}\}_{k \in \N} = \{\ve z(k,0) \}_{k \in \N}$, from Proposition \ref{problocchi} it follows that, denoting by $K_0$ a set of indices such that $\{\ve x^{(k)}\}_{k\in K_0}$ converges to $\bar{\ve x}$ and $\{\s_1^{(k,0)}\}_{k\in K_0}$ converges to some $\bar \s_1^0\in S_1$, we have $\p_1(\bar\x;\hi 1{\bar \vesi_1^0}) = \bar \x _1 $ and $\lim_{k\rightarrow \infty,k\in K_0} \|\ve d_1^{(k,0)}\|=0$.\\
From Step 2.3.4 of Algorithm \ref{AlgoBlocchi}, it follows that $\lim_{k\rightarrow \infty,k\in K_0}\|\ve x_1^{(k,1)}-\ve x_1^{(k)}\|=0$, i.e., $\ve{\bar x}_1$ is a limit point also for the sequence $\{\ve x^{(k,1)}_1\}_{k \in \N}$.\\
Introducing a subset of indices $K_1\subseteq K_0$ such that the sequence $\{\x^{(k,1)}_1\}_{k\in K_1}$ converges to $\bar \x_1$ and $\{\s_1^{(k,1)}\}_{k\in K_1}$ converges to some $\bar \s_1^1$, we have
\begin{eqnarray*}
\lim_{k\rightarrow \infty, k\in K_1} \ve d_1^{(k,1)}&=&\lim_{k\rightarrow \infty,k\in K_1}  \p_1((\ve x_1^{(k,1)},\ve x_2^{(k)},...,\ve x_m^{(k)});\hi 1{\vesi_1^{(k,1)}})-\ve x_1^{(k,1)}\\
&=&\p_1(\bar \x; \hi 1{\bar \vesi_1^1})-\bar \x_1=0,
\end{eqnarray*}
where the second equality follows from the continuity of the generalized projection operator and the third one is a consequence of Proposition \ref{problocchi}.\\
Using the same arguments, by induction on $\ell$ we can conclude that, for each $\ell =0,...,L_1-1$, there exists a suitable subset of indices $K_\ell$ such that
$\lim_{k\rightarrow \infty, k\in K_\ell} \ve d^{(k,\ell)}_1 =0$ and we obtain
\begin{equation*}
\|\ve x_1^{(k+1)}-\ve x_1^{(k)}\|\leq \sum_{\ell=0}^{L_1^{(k)}} \lambda_1^{(k,\ell)}\ve \|\d_1^{(k,\ell)}\|
\leq \sum_{\ell=0}^{L_1} \lambda_1^{(k,\ell)}\ve \|\d_1^{(k,\ell)}\|\xrightarrow{k\rightarrow \infty,k\in \bar K_1} 0,
\end{equation*}
where $\bar K_1 = \cap_{\ell = 0}^{L_1-1} K_\ell$. Thus, the point $\bar \x$ is a limit point also for the sequence $\{\ve z(k,1)\}_{k \in \N} = \{(\ve x^{(k+1)}_1,\ve x_2^{(k)},...,\ve x_m^{(k)})\}_{k \in \N}$, and Proposition \ref{problocchi} ensures that $\p_2(\bar\x;\hi 2{\bar \vesi_2^0}) = \bar \x _2 $ for some $\bar \s_2^0\in S_2$.\\
Proceeding by induction on $i$ and employing the same arguments used for $i=1$, we prove that $\ve{\bar x}$ is a limit point of the sequences $\{\ve z(k,i)\}_{k \in \N}$ for any $i=1,...,m-1$. As a result of this, invoking again Proposition \ref{problocchi}, we can conclude that for any $i=1,...,m$ there exist $\s_i\in S_i$ such that $\p_i(\bar\x;h_{\s_i}^i) = \bar \x_i $. Therefore, by Lemma \ref{lem:blocchi} (i) we can conclude that $\ve{\bar x}$ is a stationary point of problem \eqref{minf}.
\endproof
\section{Conclusions}\label{sec:concl}
In this paper we address the general problem of the constrained minimization of a differentiable function in which the unknown can be partitioned in blocks, each with a convex and closed feasible set. In order to address this problem, we considered block coordinate first order methods exploiting suitable descent directions based on very general projection operators. In particular, we introduce a class of generalized projection operators based on non Euclidean metrics, which includes as special cases Bregman projections, proximity and proximal gradient operators. Our approach combines the properties of these generalized projections with those of the Armijo linesearch strategy to obtain a generalized gradient descent method able to produce a sequence of iterates whose limit points are stationary.\\
Future work will include a generalization of these results to nonsmooth objective functions, the analysis of suitable strategies to design the parameters defining the metric functions and the extensive application of the proposed optimization approaches in real-world problems in astronomy and microscopy.

\begin{acknowledgements}
This work has been partially supported by MIUR (Italian Ministry for University and Research), under the projects FIRB - Futuro in Ricerca 2012, contract RBFR12M3AC, and PRIN 2012, contract 2012MTE38N. The Italian GNCS - INdAM (Gruppo Nazionale per il Calcolo Scientifico - Istituto Nazionale di Alta Matematica) is also acknowledged.
\end{acknowledgements}

\bibliographystyle{spmpsci}      % mathematics and physical sciences
\bibliography{cincia_biblio}   % name your BibTeX data base

%% Non-BibTeX users please use
%\begin{thebibliography}{}
%%
%% and use \bibitem to create references. Consult the Instructions
%% for authors for reference list style.
%\end{thebibliography}

\end{document}